\documentclass[12pt]{amsart}

\usepackage{amssymb}
\usepackage{amsmath}
\usepackage{amsthm}
\usepackage{amsfonts}

\usepackage[T2A]{fontenc}
\usepackage[russian, english]{babel}

\evensidemargin\oddsidemargin

\begin{document}

\baselineskip=18pt
\setcounter{page}{1}

\newtheorem{theorem}{Theorem}
\newtheorem{lemma}[theorem]{Lemma}
\newtheorem{proposition}[theorem]{Proposition}
\newtheorem{remark}[theorem]{Remark}
\newtheorem{remarks}[theorem]{Remarks}
\newtheorem{THEO}{Theorem \!\!\!\!}
\newtheorem{CORO}{Corollary\!\!}

\renewcommand{\theTHEO}{}
\renewcommand{\theCORO}{}

\newcommand{\eqnsection}{
\renewcommand{\theequation}{\thesection.\arabic{equation}}
    \makeatletter
    \csname  @addtoreset\endcsname{equation}{section}
    \makeatother}
\eqnsection

\def\e{{\mathbb E}}
\def\FF{\mathcal{F}}
\def\LL{\mathcal{L}}
\def\k{{\mathcal K}}
\def\N{{\mathbb N}}
\def\p{{\mathbb P}}
\def\ptil{{\widehat\p}}
\def\r{{\mathbb R}}
\def\z{{\mathbb Z}}
\def\SS{\mathcal{S}}

\def\lacc{\left\{}
\def\lcr{\left[}
\def\lpa{\left(}
\def\lva{\left|}
\def\racc{\right\}}
\def\rcr{\right]}
\def\rpa{\right)}
\def\rva{\right|}

\def\dt{\delta_t^\eps}
\def\eps{\varepsilon}
\def\pxy{\p_{(x,y)}}
\def\sit{\sigma_t^\eps}
\def\rt{\rho_t^\eps}
\def\Te{T_\eps}
\def\tae{\tau_\eps}
\def\Vet{{\mathcal V}^\eps_t}
\def\Un{{\bf 1}}
\def\Xde{{\tilde X}}
\def\ee{\mathrm{e}}
\def\phid{{\tilde \varphi}}

\def\d{\, \mathrm{d}}

\def\qed{\hfill$\square$}
\def\elaw{\stackrel{d}{=}}

\title[Homogeneous Brownian functionals]
      {Chung's law for homogeneous Brownian functionals}

\author[Aim\'e Lachal and Thomas Simon]{Aim\'e Lachal and Thomas Simon}

\address{Institut national des sciences appliqu\'ees de Lyon, P\^ole
  de math\'ematiques, B\^atiment L\'eonard de Vinci, 20 avenue Albert Einstein, 69621 Villeurbanne Cedex, France. {\em E-mail address}: {\tt aime.lachal@insa-lyon.fr}}

\address{Equipe d'Analyse et Probabilit\'es, Universit\'e d'Evry-Val d'Essonne, Boulevard Fran\c{c}ois Mitterrand, 91025 Evry Cedex, France. {\em E-mail address}: {\tt tsimon@univ-evry.fr}}

\keywords{First passage time, fluctuating additive functional, law of the iterated logarithm, small ball probabilities.}

\subjclass[2000]{60F99, 60G17, 60G18, 60J55, 60J65.}

\begin{abstract} Consider the first exit time $T_{a,b}$ from a finite
  interval $[-a,b]$ for an homogeneous fluctuating functional $X$ of a
  linear Brownian motion. We show the existence of a finite positive constant $\k$ such that
$$\lim_{t\to\infty}t^{-1}\log \p[ T_{ab}> t]\; =\; -\k.$$
Following Chung's original approach \cite{Ch}, we deduce a "liminf"
law of the iterated logarithm for the two-sided supremum of $X$. This
extends and gives a new point of view on a result of Khoshnevisan and Shi \cite{KS}.

\end{abstract}

\maketitle

\begin{center}
оНЯБЪЫЮЕРЯЪ еМГН нПЯХМЦЕПС Б ВЕЯРЭ ЕЦН 60-КЕРХЪ
\end{center}

\section{Introduction}
Let $\lacc B_t, \; t\ge 0\racc$ be a linear Brownian motion starting
at 0 and $X = \lacc X_t, \; t\ge 0\racc$ be the homogeneous fluctuating
additive functional defined by
$$X_t\; =\; \int_0^t V(B_s)\, \d s, \quad t \ge 0,$$
where $V(x) = x^\alpha$ if $x\ge 0$ and $V(x) = -\lambda \vert
x\vert^\alpha$ if $x\le 0,$ for some fixed $\alpha, \lambda >0.$ The
process $X$ appears in mathematical physics as the solution of a generalized
Langevin equation involving a harmonic oscillator driven by a white
noise, and we refer to \cite{La1} and the references therein for more
details on this subject. Notice that $X$ is $(1 + \alpha/2)$-self-similar,
but has no stationary increments. In the case $\alpha = \lambda =1$, it is the integrated Brownian motion:
$$X_t\; =\; \int_0^t B_s \, \d s, \quad t\ge 0,$$
and also a Gaussian process. However, in the other cases, it is not
Gaussian any longer. For every $a,b >0$ consider the bilateral exit time
$$T_{ab} \; =\;\inf\{ t >0, \; X_t\not\in(-a,b)\}.$$
As a rule, studying the law of $T_{ab}$ is a difficult issue because
$X$ alone is not Markov, so that no spectral theory is
available. We refer however to \cite{La1} and \cite{La2} for several
distributional properties of the bivariate random variable $(T_{ab},
B_{T_{ab}})$ and for the solution to the two-sided exit problem,
i.e. the computation of the probability $\p\lcr X_{T_{ab}} = a\rcr.$
In \cite{La1}, it was also shown that the variable 
$T_{ab}$ has moments of any power, and an explicit upper bound was
given on the latter - see Proposition 7.1 therein. Before this, the
upper tails of $T_{ab}$  in the case $\alpha = \lambda = 1$ had been  
precisely investigated in \cite{KS}, with an elegant argument relying
on Chung's law of the iterated logarithm. This result
was then generalized in \cite{LS} to a broad class of Gaussian and
sub-Gaussian processes, with a different method relying on wavelet
decomposition. In this paper, we aim at extending the results of
\cite{KS} to the above non-Gaussian functionals $X,$ with a more elementary proof: 

\begin{THEO} For every $a, b >0$ there exists a finite positive constant $\k$ such that
\begin{equation}
\label{Main1}
\lim_{t\to\infty}t^{-1}\log \p[ T_{ab}> t]\; =\; -\k.
\end{equation}
\end{THEO}

This exponential tail behaviour is typical for exit-times from a
finite interval for self-similar random processes. Actually, in most
examples available, it appears that the upper tails of the variable
$T_{ab}$ are those of an exponential random variable. Some comments on
this somewhat intriguing universal behaviour are given in the last
section of \cite{LS} in the case of a sub-Gaussian symmetric process
exiting a symmetric interval. See however Example 3.3 in \cite{Sam},
where the tail behaviour is shown to be subexponential. Notice also
that the upper tails of the {\em unilateral} exit time $T_{a\infty}$ of $X$
had been thoroughly studied in \cite{Iso, IK} and exhibit an entirely different, polynomial,
behaviour which again in the framework of self-similar random
processes is typical for exit-times from a semi-finite interval.
 
Taking $a=b=1$, the estimate (\ref{Main1}) entails by
self-similarity that there exists a finite positive constant $\k'$ such that
\begin{equation}
\label{Main2}
\lim_{\eps\to 0}\eps^{-2/(\alpha +2)}\log \p[ \vert\vert X\vert\vert_\infty < \eps]\; =\; -\k',
\end{equation}
where $\vert\vert .\vert\vert_\infty$ stands for the supremum norm
over $[0,1]$. This other limit theorem is known as a small ball probability
estimate, a subject which has given rise to intensive research over the
last years, with interesting connections to different questions in
analysis, probability and statistics. We refer to \cite{LiS, ThS} for recent accounts on this topic concerning both Gaussian and Non-Gaussian processes - see also Chapter 7 in \cite{Le} for an abstract Wiener setting. Originally, this
kind of estimate had been used by Chung \cite{Ch} for random walks and Brownian motion, in connection with his celebrated law of the iterated logarithm. In
\cite{KS}, Khoshnevisan \& Shi's original approach  for integrated
Brownian motion consisted in
proving first Chung's LIL and then deduce the small deviation estimate
(\ref{Main2}). In this paper, we will
follow the more standard approach viewing Chung's LIL as a consequence
of (\ref{Main2}). Introduce the notations
$$X^*_t \; =\; \sup\{\vert X_s\vert,\, s\le t\}\quad \mbox{and}\quad f(t) \; =\; (t/\log\log t)^{(\alpha+2)/2}$$   
for every $t > {\rm e},$ and set $\k_1$ for the constant appearing in (\ref{Main1}) when $a =b =1.$
\begin{CORO}[Chung's law of the iterated logarithm] One has 
$$\liminf_{t\to +\infty} \frac{X^*_t}{f(t)}\; =\; \k_1^{(\alpha +2)/2}\qquad {\rm a.s.}$$
\end{CORO}
Notice that if we introduce the family of time-stretched functionals
$$X^n_t \; =\; \frac{X_{nt}}{(n/\log\log n)^{(\alpha +2)/2}}, \quad
t\in [0,1]$$
for every $n\ge3,$ then by a straightforward monotonicity argument our
Chung's LIL is equivalent to
$$\liminf_{n\to +\infty} \vert\vert X^n\vert\vert_\infty\; =\; \k_1^{(\alpha +2)/2}\qquad {\rm a.s.}$$
From this fact and in the spirit of Wichura's functional LIL, it is an
interesting question to determine the cluster set of the family of
processes $\{X^n,\, n\ge 1\}$ for the weak topology. This was
indeed recently investigated by Lin and Zhang \cite{LZ} for $m-$fold
integrated Brownian motion, yielding Chung's LIL for these processes
as a corollary - see Theorem 1.1 and Corollary 1.1 therein. However,
in our framework the non-linearity of the kernel $x\mapsto V(x)$ and the non-Gaussianity of $X$ makes
the situation significantly more complicated in general, as it will
already appear in our proof. Setting now
$$\Xde^n_t \; =\; \frac{X_{nt}}{(n\log\log n)^{(\alpha +2)/2}}, \quad
t\in [0,1]$$
for every $n\ge3,$ our result reads
$$\liminf_{n\to +\infty} (\log\log n)^{\alpha +2}\vert\vert \Xde^n\vert\vert_\infty\; =\; \k_1^{(\alpha +2)/2}\qquad {\rm a.s.}$$
From this fact and in the spirit of Strassen's functional LIL, it is
somewhat tantalizing to determine the set of functions $f$ such that 
\begin{equation}
\label{csak}
\liminf_{n\to +\infty} (\log\log n)^{\alpha +2}\vert\vert \Xde^n
-f\vert\vert_\infty
\end{equation}
a.s. exists, as an explicit function of $f$ and $\k_1$. In the case of Brownian
motion, this (hard) problem had been initiated by Cs\'aki \cite{Cs} and De
Acosta \cite{Ac}, hinging upon shifted Brownian small balls. Of
course, before investigating (\ref{csak}) one should
first determine  the cluster set for the weak topology of the family
of processes $\{\Xde^n,\, n\ge 1\}.$ To the best of our knowledge, no results
of this kind seem to exist even for integrated Brownian motion. 

\section{Proof of the theorem}

Fix $a, b >0$ once and for all, and introduce the notation $T = T_{ab}$ for concision. For every $x, y\in\r,$ set $\pxy$ for the law of the strong Markov process $t\mapsto (B_t, X_t)$ starting at $(x,y).$ We keep the notation $\p =\p_{(0,0)}$ for brevity. Considering the function
$$\varphi(t)\; =\; \sup\lacc\pxy[T> t], \; (x,y)\in\r\,\times\,(-a,b)\racc,$$
the simple Markov property yields for every $t, s \ge 0$ 
\begin{eqnarray*}
\varphi(t+s)  & = &  \sup\lacc\p_{(x,y)}[T> s, T > t+s], \; (x,y)\in\r\,\times\,(-a,b)\racc\\
& = &   \sup\lacc\int_\r\int_a^b \p_{(x,y)}[(B_s, X_s)\in
\d u\d v,\, T> s]\p_{(u,v)}[T> t], \;
(x,y)\in\r\,\times\,(-a,b)\racc\\
& \le & \varphi(t) \times \sup\lacc\int_\r\int_a^b \p_{(x,y)}[(B_s, X_s)\in
\d u\d v,\, T> s], \; (x,y)\in\r\,\times\,(-a,b)\racc\\
& \le & \varphi(t)\varphi(s),
\end{eqnarray*}
so that the function $\psi(t) =\log\varphi(t)$ is subadditive. Hence, there exists $\k\in [0,+\infty]$ such that
$$\lim_{t\to +\infty} t^{-1}\psi(t)\; =\; \inf_{t>0}\lpa t^{-1}\psi(t)\rpa\; =\; -\k.$$
Besides from the second equality we see that $\k >0,$ since the function $\psi$ is clearly not identically zero. This entails
\begin{equation}
\label{limsup}
\limsup_{t\to\infty}t^{-1}\log \p[ T> t]\; =\; -\k\; <\; 0.
\end{equation}
The remainder of the proof will be given in two steps. First, we will
show the finiteness of $\k$, which is usually the difficult part in
small deviation problems. In the case $\alpha = \lambda = 1$, it had
been obtained in \cite{KS} through an original yet lengthy argument
relying on random normalization and Chung's LIL. Here we will provide
two proofs which are considerably simpler. The first one adapts the
elementary arguments of Lemma 1 in \cite{Ber1} to the two-dimensional Markov
process $(B,X)$, while the second one is based on the time-substitution
method which was used in \cite{Iso} for unilateral passage times - let
us stress that its main idea relying on the
a.s. continuity of the Brownian paths was also implicitly used in
\cite{KS} p. 4258 to obtain Chung's LIL.
The latter proof is slightly more involved than the former,
nevertheless it allows to bound the constant from above - see the
Remark 1 below.  

Second, we will show that the above limit in
(\ref{limsup}) is actually a true limit, which appears to be quite more
complicated. In the Gaussian case 
$\alpha = \lambda = 1$ and for a symmetric exit interval, it is an
easy consequence of  Anderson's inequality, as already noticed in
\cite{KS}. However, no isoperimetric inequalities seem available when
$X$ is not Gaussian and this argument breaks down, so that we had to
use more bare-hand estimates, following roughly the outline of Lemma 1
in \cite{Ber1}.\\

\noindent
{\bf First proof of the finiteness of the constant.} Fixing
$A < 0< B$ and $a < c < 0 < d < b,$ introduce
the functions $\phid(t)\; =\; \inf\lacc\pxy[T> t], \; (x,y)\in [A,B]\times
[c,d]\racc$ and 
$$\Phi(t)\; =\; \inf\lacc\pxy\lcr (B_t,
X_t) \in [A,B]\times [c,d], \,T> t\rcr, \; (x,y)\in [A,B]\times
[c,d]\racc, \quad t\ge 0.$$ 
For every
$(x,y)\in [A,B]\times [c,d]$ and $t, s \ge 0$ the simple Markov
property entails
\begin{eqnarray*}
\pxy\lcr T >t+s\rcr & \ge & \pxy \lcr (B_s, X_s) \in [A,B]\times [c,d],\,
\, T >t+s\rcr \\
& = & \int_A^B\!\!\int_c^d \pxy[(B_s, X_s)\in
\d u\d v,\, T> s]\times\p_{(u,v)}[T> t]\\
& \ge & \pxy\lcr (B_s, X_s) \in [A,B]\times [c,d], \,T> s\rcr
\times\phid(t)\\
& \ge & \Phi(s)\phid(t),
\end{eqnarray*}
 so that $\phid(t+s)\ge\phid(s)\Phi(t)$ for every $t,s \ge 0.$ In
 particular 
$$\varphi(n)\;\ge \; \phid(n)\; \ge\; \Phi(1)\phid(n-1)\;\ge\;\ldots\;\ge\;
\Phi(1)^n\phid(0)\; =\; \Phi(1)^n$$
for every $n\in\N,$ which entails  $t^{-1}\psi(t)\ge\log\Phi(1)$ for
every $t>0,$ since the function $t\mapsto t^{-1}\psi(t)$ is
decreasing. We finally get
$$\k \le -\log\Phi(1).$$
Now the function $(x,y,t)\mapsto \pxy\lcr (B_t,
X_t) \in [A,B]\times [c,d], \,T> t\rcr$ is continuous on the compact
$[A,B]\times[c,d]\times[0,2]$, since it satisfies the heat equation
$$\frac{1}{2}\frac{\partial^2}{\partial x^2}\; +\;
V(y)\frac{\partial}{\partial y}\; =\; \frac{\partial}{\partial t}$$
on $\r\times (-a,b)\times\r^+.$ In particular the function 
$$(x,y)\mapsto \pxy\lcr (B_1,
X_1) \in [A,B]\times [c,d], \,T> 1\rcr$$ 
is continuous on the compact $[A,B]\times[c,d]$ and since it is obviously everywhere positive, one
has $\Phi(1) >0,$ which completes the proof.

\qed

\noindent
{\bf Second proof of the finiteness of the constant.} Let $L=\lacc
L(t,x), \; t\ge 0, \; x\in \r\racc$ be the local-time process associated with $B$ and 
$$\tau_t\; =\; \inf\{ u\ge 0,\; L(0,u) > t\},\quad t\ge 0$$
be the inverse local time of $B$ at zero. It follows easily from the
Markov property and a scaling argument that the process $t\mapsto
(\tau_t, X_{\tau_t})$ is a two-dimensional L\'evy process such that
$t\mapsto \tau_t$ is a $(1/2)-$stable subordinator and $Y : t\mapsto
Y_t = X_{\tau_t}$ a $1/(\alpha +2)$-stable process. Introducing
$$\Theta\; = \;\inf\{ t > 0, \, X_{\tau_t} \not\in (-a,b)\},$$
the a.s. continuity of Brownian trajectories yields the key-inequality
\begin{equation}
\label{stetig}
T \; \ge\; \tau_{\Theta-}\quad\mbox{a.s.}
\end{equation}
As in the proof of Theorem B in \cite{ThS1} we now decompose, for every $c >0,$
\begin{eqnarray*}
\p\lcr \Theta > t \rcr & \le & \p\lcr \tau_t < ct\rcr\; +\; \p\lcr \Theta > t , \, \tau_t \ge ct \rcr\\
& \le & \p\lcr \tau_1 < ct^{-1}\rcr\; +\; \p\lcr \tau_{\Theta-} \ge ct \rcr\\
& \le & \p\lcr \tau_1 < ct^{-1}\rcr\; +\; \p\lcr T \ge ct \rcr
\end{eqnarray*}
where we used the 2-self-similarity and the a.s. increasingness of $\tau$ in the second line, and (\ref{stetig}) in the third. By Proposition VIII.3 in \cite{Ber} and a scaling argument, there exists $\k_0$ finite such that
$$\lim_{t\to\infty} t^{-1}\log \p\lcr \Theta > t \rcr\; = \; -\k_0.$$
By Theorem 5.12.9 in \cite{BGT} there exists $\k_c\to +\infty$ as $c\to 0$ such that 
$$\lim_{t\to\infty} t^{-1}\log \p\lcr \tau_1 < ct^{-1} \rcr\; = \; -\k_c.$$
Taking $c$ small enough and putting everything together yields 
$$\liminf_{t\to\infty}t^{-1}\log \p[T> t]\; \ge\; -\k_0/c \;>\; -\infty,$$
which entails $\k < +\infty$ as desired. 

\qed

\begin{remark} \label{rm1}
    {\em The positivity parameter $\p[Y_1 >0]$ of the non
    completely asymmetric L\'evy $1/(\alpha +2)$-stable process $Y$
    had been computed in \cite{IK} - see Remark 4 therein. This makes
    it possible to bound from above the constant $\k_0$ explicitly:
    when $\lambda =1$ i.e. $Y$ is symmetric, this can be done in
    subordinating $Y$ to some Brownian motion - see Theorem 4 in
    \cite{BLM} or Proposition 8 in \cite{ThS0} - whereas when
    $\lambda\neq 1$, the same method works in
    subordinating $Y$ to some completely asymmetric stable process
    with infinite variation - see Exercise VIII.1 in \cite{Ber} - and
    using the explicit calculations of \cite{Ber1} in the completely
    asymmetric case. On the other hand,
    the scaling parameter of the stable subordinator $\tau$ is
    explicit, so that the constants  $\k_c$ are also explicit, again
    by Theorem 5.12.9 in \cite{BGT}. To put it in a nutshell, our 
    second proof allows to exhibit an explicit upper bound on $\k$,
    which we will however not include here for the sake of brevity. 
    Notice that in the case of
    integrated Brownian motion in a symmetric interval, a lower bound
    had been given in \cite{KS}, Remark 1.4. Recall also that in the 
    non-completely asymmetric framework, the exact computation of
    $\k_0$ is a long-standing and challenging problem - see \cite{Ber1, BLM, BKM} and the references therein.}

\end{remark}

\noindent
{\bf Proof of the existence of the constant.} Suppose first that 
$\alpha = \lambda = 1$ and $a=b$. Then by self-similarity and by linearity of the integral one has, for every $x,y\in\r$ and $t >0$
$$\pxy [T > t] \; =\; \p_{(xt^{-1/2},yt^{-3/2})} \lcr \vert\vert X\vert\vert_\infty < at^{-3/2}\rcr\; = \; \p\lcr \vert\vert X + f^{x,y,t}\vert\vert_\infty < at^{-3/2}\rcr$$
where $\vert\vert .\vert\vert_\infty$ stands for the supremum norm over $[0,1]$ and $f^{x,y,t} : u \mapsto yt^{-3/2} + uxt^{-1/2}.$ Hence, Anderson's inequality - see e.g. (7.5) in \cite{Le} - entails 
$$\pxy [T > t] \; = \; \p\lcr \vert\vert X + f^{x,y,t}\vert\vert_\infty < at^{-3/2}\rcr\; \le \; \p\lcr \vert\vert X\vert\vert_\infty < at^{-3/2}\rcr\; =\; \p[T > t],$$
so that $\varphi(t) =\p[T > t]$ for every $t >0,$ and (\ref{limsup})
is a true limit. Unfortunately, this simple Gaussian argument cannot
be used in general, and we will have to use a lenghtier yet elementary method, which will be divided into three lemmas. For every $\eps > 0,$ introduce 
$$\Te = \inf \lacc t > 0, \; X_t \notin (-a +\eps, b -\eps)\racc.$$

\begin{lemma} \label{max} There exist $c_1, c_2, K >0$ such that for
  every $\eps$ small enough and every $t$ large enough, there exist $x_t^\eps\in (-K,K)$ and $y_t^\eps\in(-a+\eps, b-\eps)$ such that
\begin{equation}
\label{l1}
\p_{(x_t^\eps,y_t^\eps)} [\Te >t] \; \ge \; c_2 e^{-\k(1 +c_1\eps)t}.
\end{equation}
\end{lemma}

\noindent
{\em Proof.} For every $t >0,$ we can choose $(x_t^\eps,
y_t^\eps)\in\r\times (-a+\eps, b-\eps)$ such that 
\begin{equation}
\label{xeye}
\p_{(x_t^\eps,y_t^\eps)} \lcr \Te >t+1\rcr\; \ge\; \frac{1}{2}
\sup\lacc\p_{(x,y)} \lcr \Te >t+1\rcr, \;
(x,y)\in\r\times(-a+\eps,b-\eps)\racc.
\end{equation}
Besides, by scaling and translation we have for every
$(x,y)\in\r\times(-a,b)$
$$\pxy[\Te > t+1]\; =\; \p_{(x_\eps, y_\eps)}[T_{ab} > t_\eps]$$
with the notations $x_\eps\, =\, x/(1 -2\eps/(a+b))^{1/(\alpha +2)},\; y_\eps\, =\, (y -(b-a)/2)/(1 -2\eps/(a+b)) + (b-a)/2,\; \mbox{and}\; t_\eps\, =\, (t+1)/(1-2\eps/(a+b))^{2/(\alpha+2)}.$ 
Hence, choosing some constant $c_1 > 0$ such that $1 + c_1\eps >
(1-2\eps/(a+b))^{-2/(\alpha +2)}$ for every $\eps$ small enough and by
the definition of $\k,$ we
get \\

$\sup\lacc \p_{(x,y)} \lcr \Te >t+1\rcr, \; (x, y)\in
\r\times(-a+\eps,b-\eps)\racc $
\begin{align*}
\;\; \;\;\;\; \;\; \;\;\;\;\;\; \;\;\;\; \;\;  \;\;\;\;\;\;\;\;\;\;\;\;\;\;\;\;  \;\; \;\;\;\;& = \;\,\sup\lacc \p_{(x,y)} \lcr T_{ab}
>t_\eps\rcr, \; (x, y)\in \r\times(-a,b)\racc \\
& \ge \;\,\sup\lacc \p_{(x,y)} \lcr T_{ab}
>(1+ c_1\eps)(t+1)\rcr, \; (x, y)\in \r\times(-a,b)\racc \\
& \ge \;\, e^{- \k(1+c_1\eps)(t+1)}
\end{align*}
for $t$ large enough, so that by (\ref{xeye}),
\begin{equation}
\label{c1c2}
\p_{(x_t^\eps,y_t^\eps)} \lcr \Te >t+1\rcr\; \ge \; c_2 e^{-
  \k(1+c_1\eps)t}
\end{equation}
for $t$ large enough with $c_2 = e^{-\k(1+c_1)}/2.$ Set now $K =
2(1\vee \lambda^{-1/\alpha}) (a+b)^{1/\alpha},$ fix $\eps >0$ and $t$
large enough. If $\vert x_t^\eps\vert < K,$ then by (\ref{c1c2})
$$\p_{(x_t^\eps,y_t^\eps)} \lcr \Te >t\rcr\; \ge \;\p_{(x_t^\eps,y_t^\eps)} \lcr \Te >t+1\rcr\; \ge \; c_2 e^{- \k(1+c_1\eps)t}$$
and (\ref{l1}) holds since necessarily $y_t^\eps\in(-a+\eps, b-\eps).$ If $x_t^\eps\ge K,$ then introducing the stopping time 
$$S = \inf \lacc s > 0, \; B_s = K/2\racc,$$
the definition of $K$ and the strong Markov property at $S$ entail 
$$\p_{(x_t^\eps,y_t^\eps)} \lcr \Te >t+1\rcr \; = \; \p_{(x_t^\eps,y_t^\eps)} \lcr S \le 1, \Te >t+1\rcr.$$
Indeed, if $S >1$ then $B_s \ge K/2$ for every $s\le 1,$ so that $X_1
> -a +\eps +(K/2)^\alpha > b-\eps$ and $\Te <1.$ Hence,
\begin{eqnarray*}
\p_{(x_t^\eps,y_t^\eps)} \lcr \Te >t+1\rcr & \le &
\e_{(x_t^\eps,y_t^\eps)} \lcr \Un_{\{ S \le 1, X_S\in(-a +\eps, b-\eps)\}}\p_{(K/2,X_S)} \lcr\Te >t\rcr\rcr\\
& \le &  \p_{(x_t^\eps,y_t^\eps)} \lcr S \le 1\rcr\sup\lacc  \p_{(K/2,y)} \lcr\Te >t\rcr, \; y\in (-a+\eps, b-\eps)\racc\\
& \le & \sup\lacc  \p_{(K/2,y)} \lcr\Te >t\rcr, \; y\in (-a+\eps, b-\eps)\racc.
\end{eqnarray*}
In particular, setting $c_2' = e^{-\k(1+c_1)}/4$ and ${\tilde
  x}_t^\eps =K/2,$ we see by (\ref{c1c2}) that there exists ${\tilde y}_t^\eps\in(-a+\eps, b-\eps)$ such that
$$\p_{({\tilde x}_t^\eps,{\tilde y}_t^\eps)} \lcr\Te >t\rcr \; \ge \; c_2' e^{-
  \k(1+c_1\eps)t}.$$ The case $x_t^\eps\le -K$ can be handled similarly,
and the proof of Lemma \ref{max} is complete.

\qed

We now need to show that the estimate (\ref{l1}) remains true in a suitable neighbourhood of $(x_t^\eps, y_t^\eps)$. Fixing $\eps >0$ and $(x_t^\eps, y_t^\eps)\in (-K,K)\times (-a+\eps, b-\eps)$ as above for $t$ large enough, introduce
$$\Vet\; =\; \lacc\begin{array}{ll} 
[x_t^\eps, x_t^\eps +1] \times [y_t^\eps-\eps/2, y_t^\eps +\eps/2] & \mbox{if $x_t^\eps \ge 0$,} \\

[x_t^\eps -1, x_t^\eps] \times [y_t^\eps
  -\eps/2, y_t^\eps +\eps/2] & \mbox{if $x_t^\eps < 0$.}
\end{array}\right.$$ 
The key-feature of this neighbourhood is that its volume does not
depend on $t$ and for this reason, the proof of the following lemma is a
bit technical:
\begin{lemma} \label{beta} There exists $c_3 > 0$ such that for
  every $\eps > 0$
$$\inf\lacc  \pxy \lcr T >t\rcr, \; (x,y)\in \Vet\racc \;
>\; c_3 e^{- \k(1+c_1\eps)t}, \qquad t \to +\infty.$$ 
\end{lemma}

\noindent
{\em Proof.} First, by translation invariance, one has
\begin{equation}
\label{TI}
\inf \lacc \p_{(x_t^\eps,y)} \lcr T >t\rcr, \; y\in  [y_t^\eps
  -\eps, y_t^\eps +\eps]\racc \;\ge\;  \p_{(x_t^\eps,y_t^\eps)} \lcr \Te
  >t\rcr\;\ge\; c_2e^{- \k(1+c_1\eps)t}
\end{equation}
as $t\to +\infty,$ where $c_2$ is the constant in (\ref{l1}). Suppose now $x_t^\eps \ge 0$ and
introduce the stopping time
$$\sit \; =\; \inf\lacc s >0, \; B_s = x_t^\eps\racc.$$ 
For every $(x,y)\in \Vet$ one gets from the Markov property
\begin{eqnarray*}
\pxy\lcr T >t\rcr & \ge & \pxy \lcr T > t > \sit \rcr \\
& = & \int_0^t \int_a^b \pxy\lcr \sit \in ds, X_{\sit} \in dv\rcr
\p_{(x_t^\eps, v)}\lcr T > t-s\rcr\\
& \ge & \int_0^t \int_a^b \pxy\lcr \sit \in ds, X_{\sit} \in dv\rcr
\p_{(x_t^\eps, v)}\lcr T > t\rcr\\
& \ge & \int_0^t \int_{y_t^\eps-\eps}^{y_t^\eps+\eps} \pxy\lcr \sit \in ds, X_{\sit} \in dv\rcr
\times \inf\lacc \p_{(x_t^\eps, z)}\lcr T > t\rcr, \; |z-y_t^\eps| \le \eps\racc\\
& \ge & c_2\pxy \lcr \sit \le t, \, \lva X_{\sit} - y_t^\eps \rva \le \eps\rcr  e^{- \k(1+c_1\eps)t},
\end{eqnarray*}  
where we used (\ref{TI}) in the last step. Hence, since $[-\eps/2,\eps/2]\subset [y_t^\eps -y -\eps, y_t^\eps -y
+\eps],$ it suffices to bound 
$$\pxy [ \sit \le t, \, \lva X_{\sit} - y_t^\eps \rva \le \eps]\; \ge\; \p_{(x,
  0)} \lcr \sit \le t, \; \lva X_{\sit}\rva \le \eps/2\rcr $$ 
from below. Now since $\alpha\ge 0,$ there exists $M > 0$ such that
\begin{equation}
\label{Ma}
\lva u+v\rva^\alpha \le M (\lva u\rva^\alpha + \lva v\rva^\alpha)
\end{equation}
for every $u,v\in\r,$ so that $\p_{(x,0)}$ a.s.
$$\lva X_{\sit}\rva\; \le \; M\sit \lpa x^\alpha + \lpa B_{\sit}^*\rpa^\alpha\rpa,$$
with the notation $B^*_t = \max\{\lva \beta_s\rva, \, s\le t\}$ for every
$t\ge 0$, where $\{\beta_s, \; s\ge 0\}$ is a Brownian motion starting
at zero. With the notations $\dt = x -x_t^\eps, \rt = \inf\{ s >0, \;
\beta_s = -\dt\}$ and $\theta_z = \inf\{ s >0, \;
\beta_s = z\}$ for every $z\in\r,$ this entails
\begin{eqnarray*}\p_{(x, 0)} \lcr \sit \le t, \; \lva X_{\sit}\rva \le
  \eps/2\rcr & \ge & \p\lcr \rt \le t, \; \rt \lpa \lpa B_{\rt}^*\rpa^\alpha
  + x^\alpha\rpa \le \eps/2M\rcr\\
& \ge & \p\lcr \rt \le t, \; \rho_t \lpa B_{\rt}^*\rpa^\alpha \le \eps/4M,
\; \rt x^\alpha \le \eps/4M\rcr\\
& \ge & \p\lcr \rt \le t\wedge (\eps/4M x^\alpha), \; B^*_{\rt}\le x\rcr\\
& \ge & \p\lcr \rt \le t\wedge (\eps/4M x^\alpha)\wedge \theta_x\rcr
\end{eqnarray*}
where in the fourth line we used the obvious fact that $\rt \le
\theta_{-x}$ a.s. By scaling and since $0\le\dt\le x,$ we know that 
$$\lpa\rt, \theta_x\rpa\;\elaw\; \lpa\dt\rpa^2\lpa \theta_{-1},
\theta_{x/\dt}\rpa\quad\mbox{and} \quad \theta_{x/\dt}\; \ge
\theta_1\quad \mbox{a.s.}$$
By Lemma \ref{max} we now that $x\le K+1$ and since $\dt\in [0,1],$ we
finally get 
\begin{eqnarray*}
\p_{(x, 0)} \lcr \sit \le t, \; \lva X_{\sit}\rva \le
  \eps/2\rcr & \ge & \p\lcr \theta_{-1} \le \frac{t\wedge (\eps/4M x^\alpha)}{\lpa\dt\rpa^2}\wedge \theta_{x/\dt}\rcr\\
 & \ge & \p\lcr \theta_{-1} \le (\eps/4M|K+1|^\alpha)\wedge\theta_1\rcr,
\end{eqnarray*}
which finishes the proof of Lemma \ref{beta} because the right-hand side does not depend on $t$. 

\qed

Our last lemma is intuitively obvious, but we will give a proof for
the sake of completeness. 

\begin{lemma}
\label{OM}
For every $\eps > 0$, there is a constant $c_\eps$ such that
$$\p[(B_1, X_1) \in\Vet, \; T > 1] \; > \; c_\eps$$
for every $t$ large enough.
\end{lemma}

\noindent
{\em Proof.} Fix $\eps >0$ and define $K$ as in Lemma \ref{max}. For every $(x,y)\in
(-K,K)\times (-a+\eps,b-\eps)$, there exists a piecewise linear
function $f^{x,y} : [0,1]\to\r$ starting at zero such that $f^{x,y}_1
= x + 1/2$ if $x\ge 0$ and $f^{x,y}_1 = x - 1/2$ if $x < 0$, $g^{x,y}_1 = y$ and $\tau^{x,y} >1,$ with the notations
$$g^{x,y}_t\; =\; \int_0^t V(f^{x,y}_s)\, \d s, \; t\ge 0, \quad \mbox{and}\quad \tau^{x,y} \; =\;
\inf\{ t > 0, \; g^{x,y}_t\not\in (-a, b)\}.$$
Besides, since from (\ref{Ma}) we know that a.s. $\vert\vert
X-g^{x, y}\vert\vert_{\infty} \, \le\, M \vert\vert B-f^{x,
  y}\vert\vert_{\infty}^\alpha$ for every $(x,y)$, by the definition of
$\Vet$ we have for every $t >0$
$$\lacc \vert\vert B-f^{x^\eps_t, y^\eps_t}\vert\vert_{\infty} <
(\eps/2M)^{1/\alpha}\racc\;\subset\; \lacc(B_1, X_1) \in\Vet, \; T > 1\racc.$$
On the one hand, by compacity, we can clearly choose the functions $f^{x,y}$ such that
$$M\; := \; \sup \lacc \int_0^1 \lpa \frac{\d f^{x,y}_s}{\d s}\rpa^2\d
s, \; (x, y) \in (-K,K)\times (-a+\eps,b-\eps)\racc\; <\; +\infty.$$ 
On the other hand, the Onsager-Machlup formula - see e.g. Theorem 7.8
in \cite{Le} - entails 
$$\p\lcr \vert \vert B-f^{x^\eps_t, y^\eps_t}\vert\vert_{\infty} <(\eps/2M)^{1/\alpha}\rcr \;\ge \; c_{\eps}' \exp\lcr - \frac{1}{2}\int_0^1  \lpa \frac{\d
  f^{x^\eps_t,y^\eps_t}_s}{\d s}\rpa^2\d s\rcr \; \ge \;c_{\eps}' e^{-M/2}$$
where $c_{\eps}' = \p\lcr\lva\lva B\rva\rva_\infty
< (\eps/2M)^{1/\alpha}\rcr$. Putting everything together and setting $c_{\eps} =
c_{\eps}' e^{-M/2}$ completes the proof of Lemma \ref{OM}.

\qed 

We can now conclude the proof of the existence of the constant. Fix $\eps > 0$, take $t >0$ large enough and suppose first that $x_t^\eps \ge 0$. By the Markov property at time 1, 
\begin{eqnarray*}
\p[T > t] & \ge & \p[(B_1, X_1) \in\Vet, \; T > t]
\\
& \ge &  \p[(B_1, X_1) \in\Vet, \; T > 1] \times \inf\lacc
\pxy [T > t -1], (x,y)\in\Vet \racc\\
& \ge &   c_\eps \inf\lacc
\pxy [T > t ], (x,y)\in\Vet \racc\\
& \ge &   c_\eps c_3 e^{-\k (1 +c_1\eps)t},
\end{eqnarray*}
where we used Lemma \ref{OM} in the third line and Lemma \ref{beta} in the fourth. The case $x_t^\eps
< 0$ being handled analogously, we finally obtain, for every $\eps >0,$
$$\liminf_{t\to +\infty} \frac{1}{t} \log \p[T>t]\; \ge\; 
-\k(1+c_1\eps),$$
which completes the proof in letting $\eps$ tend to 0.

\qed

\begin{remarks}{\em (a) By the self-similarity of $B$, one can actually extend the
    definition of the functionals $X$ to every $\alpha > -1$ with an absolute
    convergence of the integral. In the symmetric case $\lambda =1,$ it is even
    possible to extend this definition to every $\alpha \in (-3/2,
    1],$ viewing $X$ as a Cauchy principal value process:
$$X_t =\;\lim_{\eps\to 0} \int_0^t \Un_{\{\vert B_s\vert > \eps\}}\vert B_s\vert^\alpha {\rm sgn} (B_s) \, \d s\; =\;\lim_{\eps\to 0}\int_\r \Un_{\{\vert x\vert > \eps\}} \vert x\vert^\alpha {\rm sgn} (x) (L(t,x) -L(0,x)) \d x$$
where in the second equality we used the occupation formula and where the
second limit exists a.s. since the map $x \mapsto L(t,x)$ is a.s. $\eta$-H\"older for every
$\eta < 1/2$. For $\alpha =-1$ the process $X$ is then up to a
multiplicative constant the Hilbert transform of $L$ while for $\alpha <
-1,$ it can be viewed as a fractional derivative of $L$, and we refer
to the seminal paper \cite{BY} and Chapter 5 in \cite{Ber} for much more on this topic.\

Above, the subadditivity argument and the finiteness of the constant do not
rely on the specific value of $\alpha$, so that one gets with the same
notations
$$-\infty\; <\; \liminf_{t\to\infty}t^{-1}\log \p[ T_{ab}> t]\;\le\;\limsup_{t\to\infty}t^{-1}\log \p[ T_{ab}> t]\; <\; 0,$$
which is a weaker version of our main result. However, the {\em positivity}
assumption on $\alpha$ is crucial for Lemma \ref{max} which is the
key-step in our proof of the existence of the constant. We believe
that the limit in (\ref{limsup}) is also a true limit when $\alpha$ is
negative, but the proof requires probably less bare-hand arguments
than ours.

\vspace{2mm}

(b) In the case $\alpha = \lambda = 1,$ the process $(B,X)$ is a
Gaussian diffusion and in this case it is known that the function
$f_t : (x,y)\mapsto \pxy [T >t]$ is log-concave for every $t > 0$ - see
e.g. Proposition 1.3 in \cite{Ko}. Hence, in the case of a symmetric
interval, its maximum is attained in $(0,0)$ and this gives another
proof of the existence of the constant. Despite Theorem 1.2. in
\cite{Ko}, our intuition is that the function $f_t$ remains log-concave in general, but we were unable to prove
this. If this were true, the existence of the constant would follow
immediately in the case $\lambda = 1$ and for a symmetric
interval. Let us stress that the function $f_t$ already exhibits some
concavity properties in the framework of non-Gaussian symmetric stable processes \cite{BKM}.
}

\end{remarks}

\section{Proof of the corollary}
We will follow the outline of \cite{KS} sections 2.4 and 2.5, which
are themselves a variation on Chung's original argument. First, arguing with (\ref{Main2}) and the first Borel-Cantelli lemma exactly as in section 2.4 of \cite{KS}, one can show that 
\begin{equation}
\label{liminf}
\liminf_{t\to +\infty} \frac{X^*_t}{f(t)}\; \ge\; \k_1^{(\alpha
  +2)/2}\qquad {\rm a.s.}
\end{equation}
and we leave the verification to the reader (beware the minor
correction $R\to \log R$ on the last line p. 4258). Moreover, the arguments
of section 2.3 in \cite{KS} applied to our L\'evy
$(1+\alpha/2)$-stable process $Y : t\mapsto
X_{\tau_t}$ entail without major modification  
\begin{equation}
\label{LIL}
\liminf_{t\to +\infty} \frac{X^*_t}{f(t)}\; < \; \infty \qquad {\rm a.s.}
\end{equation}
By the 0-1 law, we know that the liminf on the left-hand side is a.s. deterministic, so that Chung's law holds by (\ref{liminf}) and (\ref{LIL}), with an unknown finite positive constant. Notice in passing that (\ref{liminf}) and (\ref{LIL}) give also a third proof of the finiteness of $\k$ in the symmetric case $a =b$, which is actually Khoshnevisan \& Shi's in the case of integrated Brownian motion. 

However, to prove that
\begin{equation}
\label{limif}
\liminf_{t\to +\infty} \frac{X^*_t}{f(t)}\; \le\; \k_1^{(\alpha
  +2)/2}\qquad {\rm a.s.}
\end{equation}
we will have to modify slightly the arguments of section 2.5 in \cite{KS}, since the kernel $x\mapsto V(x)$ is not linear in general. Fixing a small $\eps >0,$ introduce the numbers $t_n = n^{4n},$ $s_n = n^{4n +3}$ and $y_n = (1 +2\eps)\k_1^{(\alpha +2)/2}f(t_n)$ for every $n\ge 1$. Define the sequence of stopping times
$$S_0 \; =\; 0\quad\mbox{and}\quad S_n =\inf\{t > t_n + S_{n-1}, \; B_t = 0\}, \quad n\ge 1.$$ 
Finally, consider the events
$$E_n\; =\; \lacc \sup_{S_n\le t\le t_{n+1} + S_n}\lva\int_{S_n}^t
V(B_s) \d s\rva \, < \, y_{n+1}\racc\quad\mbox{and}\quad F_n\; =\;
\lacc S_n < s_n + S_{n-1}\racc$$
for every $n\ge 1.$ On the one hand, setting $r_n = s_n
-t_n$, $\p_x$ for the law of $B$ starting at $x$, and resuming the
notations of Lemma \ref{beta}, the strong Markov property, the symmetry of Brownian motion and a scaling argument yield
\begin{eqnarray*}
\p [F_n^c] & = & \int_\r \p \lcr B_{S_{n-1} +t_n} \in \d x\rcr\p_x\lcr
\theta_0 > r_n\rcr \\
& = & \int_\r \p \lcr B_{t_n} \in \d x\rcr\p\lcr B_t < \vert x\vert,\;\forall\, t\le r_n\rcr \\
& = & \int_\r \p \lcr B_1 \in \d u\rcr\p\lcr B_t < \vert u\vert
\sqrt{t_nr_n^{-1}},\;\forall\, t\le 1\rcr \\
& \sim & c\sqrt{t_nr_n^{-1}}\;\sim\; c n^{-3/2}, \qquad n\to\infty
\end{eqnarray*}
for some positive finite constant $c$, so that
$$\sum_{n\ge 1} \p [F_n^c]\; <\; +\infty.$$
By the Borel-Cantelli lemma, for almost every $\omega$ there exists
$n_0(\omega)$ such that 
$$S_n (\omega)\; < \; S_{n_0(\omega)}(\omega) \; +\; s_{n_0(\omega)
  +1}\; +\; \cdots\; +\; s_n$$
for every $n > n_0(\omega).$ Hence, by the definition of $s_n,$ there
exists $n_1(\omega) > n_0(\omega)$ such that
\begin{equation}
\label{omega}
S_n (\omega)\; < \; 2s_n
\end{equation}
for every $n\ge n_1(\omega).$ On the other hand, since 
$$E_n\; =\; \lacc \sup_{0\le t\le t_{n+1}}\lva\int_0^t V(B_{S_n +s} -
B_{S_n}) \d s\rva \, < \, y_{n +1}\racc,$$
it follows readily from the strong Markov property and the definition of $S_n$ that the events $E_n$ are mutually independent. Besides, using (\ref{Main2}) and reasoning exactly as in \cite{KS} p. 4259 entails 
$$\sum_{n\ge 1} \p [E_n]\; =\; +\infty.$$
By the second Borel-Cantelli lemma, an infinity of events $E_n$ occur
a.s. and by (\ref{omega}), we know that a.s. eventually
$[2s_n,t_{n+1}]\subset [S_n,t_{n+1}+S_n]$. This entails
$$\sup_{2s_n\le t\le t_{n+1}}\lva\int_{S_n}^t V(B_s) \d s\rva \; < \; (1+2\eps)\k_1^{(\alpha +2)/2}f(t_{n+1})\qquad\mbox{i.o.}$$
By Khintchine's LIL for Brownian motion, 
$$\liminf_{n\to +\infty}\frac{1}{f(t_{n+1})}\lva
\int_{S_n}^{2s_n} V(B_s) \d s\rva \; \le \;
\liminf_{n\to +\infty}\frac{2(1\vee\lambda)s_n}{f(t_{n+1})}\lpa
B^*_{s_n}\rpa^\alpha\; =\; 0\qquad\mbox{a.s.}$$
Putting everything together and letting $\eps\to 0$ yields
\begin{equation}
\label{LIL1}
\liminf_{n\to +\infty} \frac{1}{f(t_n)}\sup_{2s_{n-1}\le t\le t_n}\lva\int_{2s_{n-1}}^t V(B_s) \d s\rva \; \le \; \k_1^{(\alpha +2)/2}\qquad\mbox{a.s.}
\end{equation}
Finally, we know from (\ref{LIL}) that 
$$\frac{X^*_{2s_{n-1}}}{f(t_n)}\; \to\; 0\qquad\mbox{a.s.}$$
which together with (\ref{LIL1}), the usual monotonicity argument, and
the fact that a.s.
$$X^*_{t_n}\; \le\; X^*_{2s_{n-1}}\; +\; \sup_{2s_{n-1}\le t\le
  t_n}\lva\int_{2s_{n-1}}^t V(B_s) \d s\rva,$$ 
yields (\ref{limif}) as desired.

\qed

\end{document}